\documentclass[a4paper,10pt]{article}
\usepackage{amsfonts,amsthm,amssymb}

\usepackage[all]{xy}

\title{A universal property of the monoidal 2-category of cospans of finite linear orders and surjections}
\author{M. Menni\footnote{Funded by Conicet, Universit\'a dell'Insubria, ANPCyT and  Lifia.}  \and 
N.~Sabadini\footnote{Funded by Universit\'a dell'Insubria and the Italian Government PRIN project ART ({\em Analisi di sistemi di Riduzione mediante sistemi di Transizione})} \and R.~F.~C.~Walters\footnotemark[2]}

\newcommand{\DD}{\ensuremath{\mathcal D}} 
\newcommand{\CC}{\ensuremath{\mathcal C}} 
\newcommand{\EE}{\ensuremath{\mathcal E}} 

\newcommand{\nats}{\mathbb{N}}

\newcommand{\opCat}[1]{{#1}^{op}}

\newcommand{\psh}[1]
{\widehat{#1}}

\newcommand{\Lin}{\ensuremath{\mathbf{Lin}}}  
\newcommand{\sLin}{\ensuremath{\mathbf{sLin}}} 
\newcommand{\iLin}{\ensuremath{\mathbf{iLin}}} 

\newcommand{\tensor}{\otimes}
\newcommand{\unit}{I}
\newcommand{\ord}[1]{#1} 

\newcommand{\cospan}[1]{\ensuremath{\mathbf{cospan}(#1)}}
\newcommand{\cospanZ}[1]{\ensuremath{\mathbf{cospan}_0(#1)}}

\newcommand{\yy}{\ensuremath{\mathbf{y}}}
\newcommand{\zz}{\ensuremath{\mathbf{z}}}

\newcommand{\vComp}{\cdot}  
\newcommand{\hComp}{\ast}  

\newcommand{\catc}[1]{\ensuremath{\tau_{#1}}}

\theoremstyle{plain}
\newtheorem{theorem}{Theorem}[section]
\newtheorem{corollary}[theorem]{Corollary}
\newtheorem{proposition}[theorem]{Proposition}
\newtheorem{lemma}[theorem]{Lemma}

\theoremstyle{definition}
\newtheorem{definition}[theorem]{Definition}

\newtheorem{remark}[theorem]{Remark}

\begin{document}

\maketitle

\begin{abstract}
We prove that the monoidal 2-category of cospans of finite linear orders and surjections is the universal monoidal category with an object $X$ with a semigroup and a cosemigroup structures, where the two structures satisfy a certain 2-dimensional separable algebra condition.
\end{abstract}


\section{Introduction}

Universal properties of cospan-like categories have been studied in geometry and computer science. 
For example, the category of 2-cobordisms has been shown  to be the universal symmetric monoidal category with a symmetric Frobenius algebra (see \cite{KockJ2004} for an exposition and references). Further, Lack showed in \cite{Lack04} that the category of cospans of finite sets is the universal symmetric monoidal category with a symmetric separable algebra. Rosebrugh, Sabadini and Walters showed in \cite{RoseSabaWalters05} a similar property of the category of cospans of finite graphs.

The aim of this paper is to make a first step in extending these results to the 2-dimensional structure of cospans.
To concentrate attention we avoid symmetries and find that a very natural extension of Lack's work characterizes the 2-category of cospans of monotone surjections between totally ordered sets, in the world of not-necessarily symmetric monoidal 2-categories.

   Part of the work involves describing universal properties of bicategories of cospans. Work along these lines has been already done in \cite{Hermida2000} and \cite{DawsonParePronk04} and it is possible that some of the results in this paper can be obtained as a byproduct of the work done in the papers just mentioned. 
On the other hand, our present concern allows us to make some simplifying assumptions and we have decided to prove the universal properties we need without appealing to more general work. We hope that the more concrete proofs presented here will make our work more accesible and, at the same time, allow to see more clearly into the combinatorics of the structures involved.

Another relevant work is \cite{RW2002}, which is however concerned with categories rather than 2-categories.

\section{The universal semigroup}

For each $n$ in $\nats$, let ${\ord{n}}$ be the total order ${\{ 0, 1, \ldots, n-1\}}$. (So that $\ord{0}$ is the empty total order.) 
Denote by $\Lin$ the category whose objects are totally ordered sets ${\ord{n} = \{0 < \ldots < n-1\}}$ with $n$ in $\nats$ and whose morphisms are monotone functions between these orders.

   Ordinal addition is a functor ${+:\Lin\times \Lin \rightarrow \Lin}$ which together with the initial object $\ord{0}$
induces a strict monoidal category ${(\Lin, +, \ord{0})}$. 
This monoidal category is presented in detail in Section~{VII.5} of \cite{maclane} where, in particular, it is proved that ${(\ord{1}, \nabla: \ord{1} + \ord{1} \rightarrow \ord{1}, !:\ord{0}\rightarrow \ord{1})}$ is the {\em universal monoid} in the sense that for any strict monoidal category ${(\CC, \tensor, \unit)}$ together with a monoid  ${(C, m:C\tensor C\rightarrow C, u:\unit\rightarrow C)}$ in $\CC$ there exists a unique strict monoidal functor ${(\Lin, +, \ord{0})\rightarrow (\CC, \tensor, \unit)}$ which maps the monoid 
${(\ord{1}, \nabla, !)}$ to
${(C, m, u)}$. (See Proposition~1 loc. cit.)

   Now let $\sLin$ be the subcategory of $\Lin$ determined by the surjective maps. The monoidal structure on $\Lin$
restricts to $\sLin$ and exercise~3(b) of Section~{VII.5} of \cite{maclane} states that ${(\sLin, +, \ord{0})}$ has the following universal property. A {\em semigroup} in ${(\CC, \tensor, \unit)}$ is defined to be a pair
 ${(C, m:C\tensor C\rightarrow C)}$ such that $C$ is an object of $\CC$ and $m$ is associative. 
Then ${(1, \nabla)}$ is the universal semigroup.

   The main results of this paper will also use the following results concerning pushouts in $\sLin$ and their interaction with the tensor $+$. First let us say that a category has {\em strict pushouts} if every diagram 
${y\leftarrow x\rightarrow z}$ in the category can be completed to a unique pushout square.

\begin{lemma}\label{LemStrictPOinSLin}
$\sLin$ has strict pushouts.
\end{lemma}
\begin{proof}
It is straightforward to see that the pushout (in the category of finite ordinals and all funtions) of two surjections ${p\leftarrow m \rightarrow n}$ yields a pushout ${p\rightarrow q \leftarrow n}$ in the category of ordinals and surjective functions. Among these pushouts, there is a unique one making the function ${m\rightarrow q}$ order preserving. 
\end{proof}

   Also, let us say that {\em pushouts interact with $\tensor$} in ${(\CC, \tensor, \unit)}$ if it holds that whenever the left and middle squares below are pushouts then so is the one on the right.
$$\xymatrix{
x \ar[d]_-{f} \ar[r]^-{g} & z \ar[d]^-{p_1} && x \ar[d]_-{f'} \ar[r]^-{g'} & z \ar[d]^-{p'_1} &&
   x \tensor x' \ar[d]_-{f\tensor f'}\ar[r]^-{g \tensor g'} & z \tensor z' \ar[d]^{p_1 \tensor p'_1} \\
y \ar[r]_-{p_0}           & P               && y \ar[r]_-{p'_0}           & P'                &&
   y \tensor y' \ar[r]_-{p_0 \tensor p'_0} & P \tensor P'
}$$

\begin{lemma}\label{LemPOinteract} Pushouts interact with $+$ in ${(\sLin, +, \ord{0})}$.
\end{lemma}
\begin{proof}
Obvious.
\end{proof}

\section{Cospans}

   In this section let ${\CC}$ be a category with strict pushouts. Then $\cospan{\CC}$ has the structure of a 2-category
and there are obvious functors ${\yy:\CC\rightarrow \cospan{\CC}}$ and ${\zz:\opCat{\CC}\rightarrow \cospan{\CC}}$ such that for every $C$ in $\CC$, ${\yy C = \zz C}$. For every arrow ${f:C\rightarrow C'}$ in $\CC$, ${\yy f}$ is the cospan
${(f:C\rightarrow C' \leftarrow C': id)}$ and ${\zz f}$ is the cospan ${(id:C'\rightarrow C'\leftarrow C:f)}$.

   We write composition in `Pascal' notation. So, for example, the commutative square above translates to the equation  ${(F_1 \alpha);(F_0\beta) = (F_0 p_0); (F_1 p_1)}$.

   Now let ${\DD}$ be a  2-category. Each  2-functor ${\cospan{\CC}\rightarrow \DD}$ induces by composition functors ${\CC\rightarrow \DD}$ and ${\opCat{\CC}\rightarrow \DD}$ which coincide at objects. 
In this section we describe what else is needed go the other way around.

\begin{definition}\label{DefCompatible}
A pair of functors ${F_0:\CC\rightarrow \DD}$ and ${F_1:\opCat{\CC}\rightarrow \DD}$ are called {\em compatible} if 
\begin{enumerate}
\item they coincide at the level of objects (and so we write ${F_0 X = F X = F_1 X}$)  
\item for every pushout square as on the left below,
$$\xymatrix{
X \ar[d]_-{\alpha} \ar[r]^-{\beta} & B \ar[d]^-{p_1} && F X \ar[r]^-{F_0 \beta}                    & F B \\
A \ar[r]_-{p_0}                   & P               && F A \ar[u]^-{F_1 \alpha} \ar[r]_-{F_0 p_0} & F P \ar[u]_-{F_1 p_1}
}$$
the square on the right above commutes. 
\end{enumerate}
\end{definition}

\begin{lemma}\label{LemYYZZ}
The functors $\yy$ and $\zz$ are compatible.
\end{lemma}
\begin{proof}
Straightforward.
\end{proof}

   Another simple but important fact is the following.

\begin{lemma}\label{LemTrivialCompatible}
Let ${F_0:\CC\rightarrow \DD}$ and ${F_1:\opCat{\CC}\rightarrow \DD}$ be compatible functors and let 
${G:\DD\rightarrow \EE}$ a functor. Then ${F_0 ; G}$ and ${F_1 ; G}$ are also compatible.
\end{lemma}

   Notice that in Definition~\ref{DefCompatible} we are not requiring $\DD$ to be a 2-category. For the next result let $\cospanZ{\CC}$ denote the underlying ordinary category of the 2-category ${\cospan{\CC}}$.

\begin{lemma}\label{LemOneDimExtension}
Let $\DD$ be a category and let ${F_0:\CC\rightarrow \DD}$ and ${F_1:\opCat{\CC}\rightarrow \DD}$ be functors.
Then there exists a unique ${F:\cospanZ{\CC}\rightarrow \DD}$ such that ${\yy ; F = F_0}$ and ${\zz ; F = F_1}$ if and only if $F_0$ and $F_1$ are compatible.
\end{lemma}
\begin{proof}
   One direction is trivial by  Lemmas~\ref{LemYYZZ} and \ref{LemTrivialCompatible}.
   On the other hand, assume that $F_0$ and $F_1$ are compatible. The conditions ${\yy ; F = F_0}$ and ${\zz ; F = F_1}$
determine the definition of $F$ on objects. To deal with 1-cells notice that every cospan ${p = (p_0:A\rightarrow P \leftarrow B:p_1)}$ is the result of the composition ${(\yy p_0);(\zz p_1)}$. So as $F$ must preserve composition 
${Fp = (F_0 p_0) ; (F_1 p_1): FA \rightarrow FB}$. So the definition of $F$ is forced and it remains to check that so defined $F$ is indeed a functor. Identities are preserved because $F_0$ and $F_1$ preserve them. Concerning composition, 
let ${f = (f_0:X\rightarrow A \leftarrow Y:f_1)}$ and ${g = (g_0:Y \rightarrow B \leftarrow Z:g_1)}$ be a pair of composable cospans. If we let the following square be the pushout of $f_1$ and $g_0$
$$\xymatrix{
Y\ar[d]_-{f_1} \ar[r]^-{g_0} & B \ar[d]^-{p_1} \\
A \ar[r]_-{p_0}              & P
}$$
then ${f;g}$ is the cospan ${(f_0 ; p_0) : X\rightarrow P \leftarrow Z: (g_1 ; p_1)}$.
Now calculate using compatibility (and recall that $F_1$ is contravariant):
\[ F (f;g) = F_0 (f_0; p_0) ; F_1 (g_1; p_1) = (F_0 f_0); (F_0 p_0) ; (F_1 p_1) ; (F_1 g_1) = \]
\[ = (F_0 f_0); (F_1 f_1) ; (F_0 g_0) ; (F_1 g_1) = (F f) ; (F g) \]
so the result is proved.
\end{proof}

   At greater generality, an analogous result is dealt with in Example~{5.3} of \cite{Lack04}. We prefer to be somewhat more explicit as it will allow us to see more clearly into how to extend the results one dimension up.

\section{The extension to 2-cells}

   When considering 2-categories, ${\iota_{(\_)}}$ denotes the operation providing identities for horizonal and vertical composition. That is, 2-cells of the form ${\iota_f}$ act as units for vertical composition and those of the form
${\iota_{id_A}}$ act as units for horizontal composition.
Also, vertical composition of 2-cells is denoted by $\vComp$ and horizontal one by $\hComp$. In all cases we write composition in `Pascal' order.


\begin{definition}\label{DefCompatibleFamOf2cells}
Let $\DD$ be a  2-category and ${F_0:\CC\rightarrow \DD}$ and ${F_1:\opCat{\CC}\rightarrow\DD}$ be functors. 
A {\em compatible selection of 2-cells} is a function ${\catc{(\_)}}$ that assigns to each map ${f:X\rightarrow Y}$ in $\CC$ a two cell ${\catc{f}:id_{FX} \Rightarrow (F_0 f);(F_1 f)}$ such that:
\begin{enumerate}
\item ${\catc{id_X} = \iota_{(id_{FX})}}$
\item ${\catc{\alpha ; \beta} = \catc{\alpha} \vComp ((F_0 \alpha)\hComp \catc{\beta} \hComp (F_1 \alpha))}$
\item for every pushout in $\CC$ as below
$$\xymatrix{
X \ar[d]_-{\alpha} \ar[r]^-{\beta} & B \ar[d]^-{p_1} \\
A \ar[r]_-{p_0}                   & P               
}$$
the following identities hold: 
\[ \catc{\alpha} \hComp (F_0 \beta) = (F_0 \beta) \hComp \catc{p_1} \quad\quad\quad  
(F_1 \alpha) \hComp \catc{\beta} = \catc{p_0} \hComp (F_1 \alpha). \]
\end{enumerate}
\end{definition}

   The idea is, of course, that a compatible selection of 2-cells is exactly what is needed to extend
Lemma~\ref{LemOneDimExtension} to two dimensions. But before we prove the result let us prove a couple of technical lemmas.

   First notice that each ${\alpha:A\rightarrow B}$ in $\CC$ induces a 2-cell 
${\overline{\alpha}:id_A\Rightarrow (\yy \alpha);(\zz \alpha)}$.

\begin{lemma}\label{LemBasicAsignment}
The assignment ${\alpha\mapsto \overline{\alpha}}$ is compatible with $\yy$ and $\zz$.
\end{lemma}
\begin{proof}
Straightforward.
\end{proof}

   Now we need a result analogous to Lemma~\ref{LemTrivialCompatible}.

\begin{lemma}\label{Lem2FuncInduceComposite2Cells}
Let $\CC$ be a category and let $\DD$ and $\EE$ be 2-categories. Moreover, let ${F_0:\CC\rightarrow\DD}$ and 
${F_1:\opCat{\CC}\rightarrow\DD}$ be compatible functors and let ${G:\DD\rightarrow \EE}$ be a  2-functor. 
If $\catc{(\_)}$ is a selection of 2-cells compatible with $F_0$ and $F_1$ then ${G\catc{(\_)}}$ is compatible with 
${F_0 ; G}$ and ${F_0 ; G}$.
\end{lemma}
\begin{proof}
Easy.
\end{proof}

\begin{lemma}\label{LemHorComp}
Let ${p_0}$ and ${p_1}$ be the pushout of $\alpha$ and $\beta$ as in Definition~\ref{DefCompatibleFamOf2cells} and 
let ${\gamma = \alpha ; p_0 = \beta ; p_1}$. Then ${\catc{\alpha} \hComp \catc{\beta} = \catc{\gamma}}$.
\end{lemma}
\begin{proof}
Calculate:
\[ \catc{\alpha ; p_0} = \catc{\alpha} \vComp ((F_0 \alpha)\hComp \catc{p_0} \hComp (F_1 \alpha)) = 
\catc{\alpha} \vComp ((F_0 \alpha)\hComp (F_1 \alpha) \hComp \catc{\beta} ) = \]
\[ = (\catc{\alpha} \hComp \iota_{X}) \vComp ( [(F_0 \alpha)\hComp \iota_{A} \hComp (F_1 \alpha)] \hComp \catc{\beta} ) =
(\catc{\alpha} \vComp [(F_0 \alpha)\hComp \catc{id_A} \hComp (F_1 \alpha)]) \hComp (\iota_X \vComp \catc{\beta}) = \]
\[ = \catc{(\alpha; id_A)} \hComp \catc{\beta} = \catc{\alpha} \hComp \catc{\beta} \]
\end{proof}

We can now prove the following.

\begin{proposition}\label{PropTwoDimExtension}
Let $\DD$ be a  2-category, let ${F_0:\CC\rightarrow \DD}$ and ${F_1:\opCat{\CC}\rightarrow \DD}$ be functors and let
${\catc{(\_)}}$ be a function assigning a 2-cell to each map in \CC. Then the following are equivalent:
\begin{enumerate}
\item there exists a unique   2-functor ${F:\cospan{\CC}\rightarrow \DD}$ such that the equations 
${\yy ; F = F_0}$, ${\zz ; F = F_1}$ and ${\catc{(\_)} = F\overline{(\_)}}$ hold; (here and for the rest of the paper 
${F\overline{(\_)}}$ is denoting the operation that to each 1-cell $f$ in $\CC$ assigns the 2-cell ${F(\overline{f})}$)
\item $F_0$ and $F_1$ are compatible and $\catc{(\_)}$ is a selection of 2-cells that is compatible with them.
\end{enumerate}
\end{proposition}
\begin{proof}
Assume that the first item holds. Lemma~\ref{LemOneDimExtension} implies that ${\yy ; F = F_0}$ and ${\zz ; F = F_1}$ are compatible. Lemmas~\ref{LemBasicAsignment} and \ref{Lem2FuncInduceComposite2Cells} show that $\catc{(\_)}$ is a selection of 2-cells that is compatible with them.

To prove the converse notice that we can apply Lemma~\ref{LemOneDimExtension} again to conclude that there exists a unique ordinary functor ${F:\cospanZ{\CC}\rightarrow \DD}$ such that ${\yy ; F = F_0}$ and ${\zz ; F = F_1}$.
So we are left to show that this $F$ extends uniquely to a 2-functor in such a way that 
${\catc{(\_)} = F\overline{(\_)}}$ holds.

First assume that the functor $F$ does extend to a   2-functor and consider an arbitrary 2-cell $\alpha$ as below.
$$\xymatrix{
                                & A \ar[dd]^-{\alpha} & \\
X \ar[ru]^-{f_0} \ar[rd]_-{g_0} &                     & \ar[lu]_-{f_1} \ar[ld]^-{g_1} Y \\
                                & B  
}$$
If we denote the cospans ${(f_0:X\rightarrow A \leftarrow Y:f_1)}$ and ${(g_0:X\rightarrow A \leftarrow Y:g_1)}$ by ${f:X\rightarrow Y}$ and $g$ respectively then it is easy to see that, for $\alpha$ considered as a 2-cell 
${f\Rightarrow g}$, ${\alpha = (\yy f_0) \hComp \overline{\alpha} \hComp (\zz f_1)}$.
So ${F(\alpha:f\Rightarrow g) = F (\yy f_0) \hComp F \overline{\alpha} \hComp F(\zz f_1)}$ and hence the definition of $F$ at the level of 2-cells is completely determined by 
${F(\alpha:f\Rightarrow g) = (F_0 f_0) \hComp \catc{\alpha} \hComp (F_1 f_1)}$.

   Finally we are left to prove that if we define ${F(\alpha:f\Rightarrow g)}$ to be the 2-cell 
${(F_0 f_0) \hComp \catc{\alpha} \hComp (F_1 f_1)}$ as above then we do obtain a   2-functor.
\begin{enumerate}
\item ${F \iota_f = \iota_{Ff}}$
\[ F \iota_f = (F_0 f_0) \hComp \catc{id_A} \hComp (F_1 f_1) = (F_0 f_0) \hComp \iota_{id_A} \hComp (F_1 f_1) = \]
\[ = (F_0 f_0) \hComp (F_1 f_1) = \iota_{(F_0 f_0);(F_1 f_1)} = \iota_{(Ff)} \]

\item Consider  maps ${\alpha:A\rightarrow B}$ and 
${\beta:B\rightarrow C}$ inducing 2-cells in unique posible way (starting from $f$). Then calculate
\[ F(\alpha \vComp \beta) = (F_0 f_0) \hComp \catc{(\alpha; \beta)} \hComp (F_1 f_1) = \]
\[ = (F_0 f_0) \hComp (\catc{\alpha} \vComp ((F_0 \alpha) \hComp \catc{\beta} \hComp (F_1 \alpha))) \hComp (F_1 f_1) = \]
\[ = (\iota_{(F_0 f_0)} \vComp \iota_{(F_0 f_0)}) \hComp (\catc{\alpha} \vComp ((F_0 \alpha) \hComp \catc{\beta} \hComp (F_1 \alpha))) \hComp (F_1 f_1) = \]
\[ = [(\iota_{(F_0 f_0)} \hComp \catc{\alpha}) \vComp (\iota_{(F_0 f_0)}\hComp ((F_0 \alpha) \hComp \catc{\beta} \hComp (F_1 \alpha))) ] \hComp (F_1 f_1) = \]
\[ = [(F_0 f_0) \hComp \catc{\alpha}) \vComp ((F_0 f_0)\hComp (F_0 \alpha) \hComp \catc{\beta} \hComp (F_1 \alpha)) ] \hComp 
(\iota_{(F_1 f_1)} \vComp \iota_{(F_1 f_1)}) = \]
\[ = [(F_0 f_0) \hComp \catc{\alpha} \hComp (F_1 f_1)] \vComp
[(F_0 f_0)\hComp (F_0 \alpha) \hComp \catc{\beta} \hComp (F_1 \alpha)) \hComp (F_1 f_1)]  = (F\alpha) \vComp (F\beta) \]

\item Preservation of horizontal composition. Suppose we have 2-cells ${\alpha:f\Rightarrow f'}$ and 
${\beta:g\Rightarrow g'}$ as in the diagram below.

$$
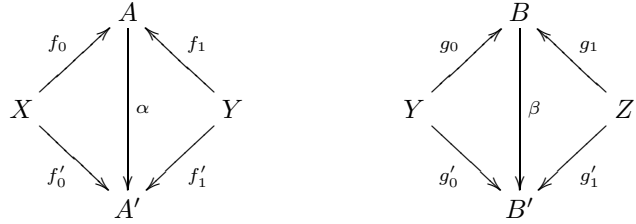
\xymatrix{
                                & A \ar[dd]^-{\alpha} &                                      && 
  & B \ar[dd]^-{\beta} & \\
X \ar[ru]^-{f_0} \ar[rd]_-{f'_0} &                     & \ar[lu]_-{f_1} \ar[ld]^-{f'_1} Y && 
  Y \ar[ru]^-{g_0} \ar[rd]_-{g'_0} &                     & \ar[lu]_-{g_1} \ar[ld]^-{g'_1} Z \\
                                & A'                  &                                 &&
                                  & B'                  &                                
}$$
In order to calculate ${\alpha\hComp \beta}$ calculate the following pushout and resulting map (every small square is a push out).
$$
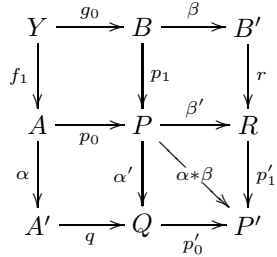
\xymatrix{
Y \ar[d]_-{f_1} \ar[r]^-{g_0}    & B \ar[d]^-{p_1} \ar[r]^-{\beta}                                 & B' \ar[d]^-{r}   \\
A \ar[d]_-{\alpha} \ar[r]_-{p_0} & P \ar[rd]|-{\alpha\hComp\beta} \ar[d]_-{\alpha'} \ar[r]^-{\beta'} & R \ar[d]^-{p'_1} \\
A' \ar[r]_q                      & Q \ar[r]_-{p'_0}                                                & P'
}$$
Now use Lemma~\ref{LemHorComp} to calculate:
\[ F(\alpha\hComp\beta) = F_0(f_0 ; p_0) \hComp \catc{\alpha\hComp\beta} \hComp F_1(g_1 ; p_1) = \]
\[ = (F_0 f_0) \hComp (F_0 p_0) \hComp \catc{\alpha'} \hComp \catc{\beta'} \hComp (F_1 p_1) \hComp (F_1 g_1) = \]
\[ = (F_0 f_0) \hComp \catc{\alpha} \hComp (F_0 p_0) \hComp (F_1 p_1) \hComp \catc{\beta} \hComp (F_1 g_1) = \]
\[ = (F_0 f_0) \hComp \catc{\alpha} \hComp (F_1 f_1) \hComp (F_0 g_0) \hComp \catc{\beta} \hComp (F_1 g_1) = 
(F\alpha) \hComp (F\beta) \]
\end{enumerate}

\end{proof}

\section{Adding the monoidal structure}

In this section let ${(\CC,\oplus,0)}$ be a strict monoidal category with strict pushouts.
We have already seen that ${\cospan{\CC}}$ is a  2-category. We want to `extend' the tensor $\oplus$ on $\CC$ to one on $\cospan{\CC}$. More precisely we will construct a   2-functor 
${\oplus: \cospan{\CC}\times\cospan{\CC} \rightarrow \cospan{\CC}}$.

\begin{lemma}\label{LemProdsAndCospans}
There exists a 2-iso ${\cospan{\CC}\times\cospan{\CC} \rightarrow \cospan{\CC\times\CC}}$ such that the following diagram commutes.
$$\xymatrix{
\CC\times\CC \ar[rd]_-{\yy} \ar[r]^-{\yy\times\yy} & \cospan{\CC}\times\cospan{\CC}  \ar[d]^-{\cong} &
  \ar[l]_-{\zz\times\zz} \opCat{\CC}\times \opCat{\CC} \ar[d]^-{\cong}        \\
 &  \cospan{\CC\times\CC}    & \ar[l]_-{\zz} \opCat{(\CC\times \CC)}
}$$
\end{lemma}
\begin{proof}
The obvious one.
\end{proof}

   So we need only build a   2-functor ${\oplus: \cospan{\CC\times\CC} \rightarrow \cospan{\CC}}$ with the right properties.

\begin{lemma}\label{LemTensorYYZZareCompatible}
The functors 
$$\xymatrix{
\CC\times \CC \ar[r]^-{\oplus} & \CC \ar[r]^-{\yy} & \cospan{\CC} &  
\opCat{(\CC\times \CC)} \ar[r]^-{\opCat{\oplus}} & \opCat{\CC} \ar[r]^-{\zz} & \cospan{\CC}
}$$
are compatible if and only if pushouts interact with $\oplus$ in ${(\CC, \oplus, 0)}$. Moreover, in this case ${(\alpha,\beta) \mapsto \overline{\alpha \oplus \beta}}$ is a compatible selection of 2-cells.
\end{lemma}
\begin{proof}
The functors coincide at the level of objects.
Now, a pushout in ${\CC\times\CC}$ is a pair of pushouts ${\alpha; p_0 = \beta ; p_1}$ and 
${\alpha'; p'_0 = \beta' ; p'_1}$ in $\CC$. The compatibility condition reduces, in this case, to 
${(\alpha\oplus\alpha');(p_0\oplus p'_0) = (\beta\oplus\beta');(p_1\oplus p'_1)}$ being a pushout.
So the first part of the result follows.

   For the second part denote let ${\sigma_{(\alpha,\beta)} = \overline{\alpha \oplus \beta}}$ and recall that $\overline{(\_)}$ is a compatible selection of 2-cells (Lemma~\ref{LemBasicAsignment}).
It is easy to show that ${\sigma_{(id_X, id_Y)} = \iota_{id_{X\oplus Y}}}$. 
In order to check the second conditon calculate:
\[ \sigma_{(\alpha, \alpha') ; (\beta, \beta')} = \sigma_{((\alpha;\beta) , (\alpha';\beta'))} = 
\overline{(\alpha;\beta)\oplus (\alpha';\beta')} = \overline{(\alpha\oplus\alpha') ; (\beta\oplus\beta')} = \]
\[ = \overline{(\alpha\oplus\alpha')} \vComp
(\yy (\alpha\oplus\alpha') \hComp \overline{(\beta\oplus\beta')} \hComp \zz(\alpha\oplus\alpha')) =
\sigma_{(\alpha,\alpha')} \vComp (\yy (\alpha\oplus\alpha') \hComp \sigma_{(\beta, \beta')} \hComp \zz(\alpha\oplus\alpha'))\]

In order to check the final condition assume that we have a pushout in $\CC\times\CC$ as on the left below
$$\xymatrix{
\ar[d]_-{(\alpha, \alpha')} \ar[r]^-{(\beta, \beta')} & \ar[d]^-{(p_1, p'_1)} & & 
  \ar[d]_-{\alpha\oplus \alpha'} \ar[r]^-{\beta\oplus \beta'} & \ar[d]^-{p_1\oplus p'_1} \\
\ar[r]_{(p_0, p_1)} &                                                         & & \ar[r]_{p_0\oplus p_1} &
}$$
then the square on the right above is a pushout because of interaction. Then calculate:
\[ \sigma_{(\alpha, \alpha')} \hComp \yy (\beta\oplus\beta') = \overline{\alpha\oplus\alpha'} \hComp \yy (\beta\oplus\beta') =
\yy (\beta\oplus\beta') \hComp \overline{p_1\oplus p'_1} = \yy (\beta\oplus\beta') \hComp \sigma_{(p_1, p'_1)} \]

The equation ${\zz(\alpha\oplus \alpha') \hComp \sigma_{(\beta, \beta')} = \sigma_{(p_0, p'_0)} \zz(\alpha\oplus \alpha')}$ is dealt with similarly.
\end{proof}

\begin{proposition}\label{PropMonoidal2cat}
Let ${(\CC,\oplus,0)}$ be a strict monoidal category with strict pushouts then the following are equivalent:
\begin{enumerate}
\item there exists a unique   2-functor 
\[ \oplus:\cospan{\CC}\times\cospan{\CC}\rightarrow \cospan{\CC} \]
such that that the following diagrams commute
$$\xymatrix{
\CC\times\CC \ar[d]_-{\yy\times\yy} \ar[r]^-{\oplus} & \CC \ar[d]^-{\yy}        \\
\cospan{\CC}\times\cospan{\CC} \ar[r]_-{\oplus}      & \cospan{\CC} 
}$$
$$\xymatrix{
\opCat{\CC}\times\opCat{\CC} \ar[d]_-{\zz\times\zz} \ar[r]^-{\opCat{\oplus}} & \opCat{\CC} \ar[d]^-{\zz} \\
\cospan{\CC}\times\cospan{\CC} \ar[r]_-{\oplus}                      & \cospan{\CC}
}$$
and such that ${\overline{\alpha} \oplus \overline{\beta} = \overline{\alpha \oplus\beta}}$.

\item Pushouts and $\oplus$ interact in ${(\CC,\oplus, 0)}$.
\end{enumerate}
%
In  this  case, the resulting 
structure ${(\cospan{\CC}, \oplus, 0)}$ is a  monoidal 2-category and the 
functors ${\yy}$ and $\zz$ extend to strict monoidal ${(\CC, \oplus, 0)\rightarrow \cospan{\CC}}$ and
${(\opCat{\CC}, \oplus, 0)\rightarrow \cospan{\CC}}$ respectively.
\end{proposition}
\begin{proof}
Lemma~\ref{LemTensorYYZZareCompatible} together with  Proposition~\ref{PropTwoDimExtension} show that the interaction of pushouts with $\oplus$ in the structure ${(\CC,\oplus, 0)}$ is equivalent to the existence of a  2-functor 
${\oplus:\cospan{\CC\times\CC} \rightarrow\cospan{\CC}}$ satisfying a number of properties which, after precomposing with the isomorphism of Lemma~\ref{LemProdsAndCospans}, turn out to be exactly the ones in the statement of the present result.

The rest of the statement is trivial by strictness. 
\end{proof}

 In particular:

\begin{corollary}\label{CorSLinIsMon2cat}
${(\cospan{\sLin}, +, \ord{0})}$ is a strict monoidal 2-category.
\end{corollary}

\subsection{Monoidal 2-functors from $\cospan{\CC}$}

   In this section let ${(\CC, \oplus, 0)}$ be a strict monoidal category such that $\CC$ has strict pushouts that interact with $\oplus$. By Proposition~\ref{PropMonoidal2cat} we have the strict monoidal 2-category 
${(\cospan{\CC}, \oplus, 0)}$.

\begin{lemma}\label{LemMonoidalInduced}
Let ${(\DD, \tensor, \unit)}$ be a strict monoidal 2-category. Let  ${F_0:\CC\rightarrow\DD}$ and
${F_1:\opCat{\CC}\rightarrow\DD}$ be compatible functors. Finally let $\catc{(\_)}$ be a compatible selection of 2-cells.
Then the induced   2-functor ${F:\cospan{\CC}\rightarrow \DD}$ is strict monoidal 
${(\cospan{\CC}, \oplus, 0)\rightarrow (\DD, \tensor, \unit)}$ if and only if ${F_0}$ and ${F_1}$ are strict monoidal and
${\catc{\alpha\oplus\beta} =\catc{\alpha}\tensor \catc{\beta}}$.
\end{lemma}
\begin{proof}
Assume that $F$ is a strict monoidal 2-functor. Then clearly ${F_0 = \yy; F}$ is a strict monoidal functor and similarly for $F_1$. To check the condition on the selection of 2-cells just calculate:
\[ \catc{\alpha\oplus\beta} = F(\overline{\alpha\oplus\beta}) = F(\overline{\alpha}\oplus\overline{\beta}) =
F\overline{\alpha} \tensor F\overline{\beta} = \catc{\alpha} \tensor \catc{\beta} \]

Conversely, assume that the conditions stated for $F_0$, $F_1$ and $\catc{(\_)}$ hold. 
Clearly ${F(X\oplus Y) = FX \tensor FY}$ because $F$ and $F_0$ coincide at the level of objects.
Now let ${f = (f_0:X \rightarrow A \leftarrow Y:f_1)}$ and ${f' = (f'_0:X' \rightarrow A' \leftarrow Y':f'_1)}$ be 1-cells. 
Then calculate
\[ F( f \oplus f') = F(f_0\oplus f'_0:X\oplus X' \rightarrow A\oplus A' \leftarrow Y\oplus Y':f_1\oplus f'_1) = \]
\[ = F_0(f_0 \oplus f'_0) ; F_1(f_1 \oplus f'_1) = 
 ((F_0 f_0)\otimes (F_0 f'_0)) ; ((F_1 f_1)\otimes (F_1 f'_1)) = \]
\[ = ((F_0 f_0) ; (F_1 f_1)) \otimes ((F_0 f'_0) ; (F_1 f'_1)) = (Ff) \otimes (Ff') \]

Finally, consider a 2-cells $\alpha$ from $f$ and $\beta$ from $f'$ and calculate using that $\tensor$ is a   2-functor:
\[ F(\alpha\oplus\beta) = F_0(f_0\oplus f'_0) \hComp \catc{\alpha\oplus\beta} \hComp F_1(f_1\oplus f'_1) = \]
\[= ((F_0 f_0)\otimes (F_0 f'_0))\hComp (\catc{\alpha} \tensor \catc{\beta})\hComp ((F_1 f_1) \otimes (F_1 f'_1)) = \]
\[ = ( (F_0 f_0) \hComp \catc{\alpha} \hComp (F_1 f_1)) \tensor ( (F_0 f'_0) \hComp \catc{\beta} \hComp (F_1 f'_1)) = 
(F\alpha) \tensor (F\beta) \]
\end{proof}

\section{Separable semi-algebras}

In this section we introduce the fundamental 1-dimensional structure to be studied in the paper.

\begin{definition}
Let ${(\DD, \tensor, \unit)}$ be a strict monoidal category. A {\em bi-semigroup} ${(X, \nabla, \Delta)}$ is an object  $X$ in $\DD$ together with morphisms ${\nabla:X\tensor X \rightarrow X}$ and ${\Delta:X\rightarrow X\tensor X}$ such that 
${(X, \nabla)}$ is a semigroup and ${(X, \Delta)}$ is a `co-semigroup' in the sense that $\Delta$ is coassociative.
\end{definition}

It is useful to have a graphical notation for expressions involving $\nabla$ and $\Delta$.
A couple of examples will suffice to introduce it. Consider a bi-semigroup ${(X, \nabla, \Delta)}$. The identity on $X$ will be denoted by a straight line. On the other hand, ${id_X \tensor id_X}$ will be denoted by two parallel horizontal lines. More importantly, $\nabla$ will be denoted as in the left diagram below
$$
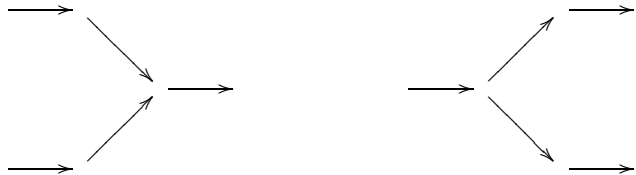
\xymatrix{
\ar[r] & \ar[rd] &        & &&        &                & \ar[r] & \\
       &         & \ar[r] & && \ar[r] & \ar[ru]\ar[rd] &        & \\ 
\ar[r] & \ar[ru] &        & &&        &                & \ar[r] & \\
}$$
and $\Delta$  will be denoted as on the right above. So that, for example, the diagram 
$$
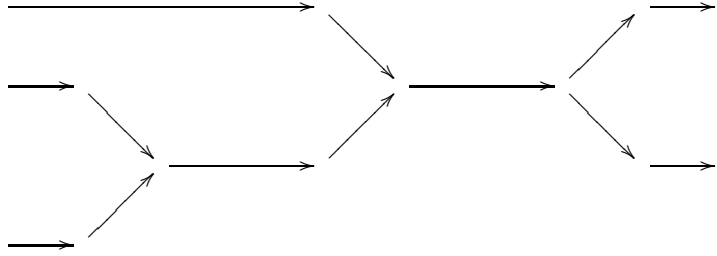
\xymatrix{
\ar[rrrr] &       &        &         & \ar[rd] &        &         &                & \ar[r] & \\
\ar[r] & \ar[rd] &        &        &         & \ar[rr] &  & \ar[ru]\ar[rd] & \\ 
       &         & \ar[rr] &       & \ar[ru] &        &         &                & \ar[r] & \\
\ar[r] & \ar[ru] &  \\
}$$
represents the expression ${(id\tensor\nabla);\nabla;\Delta}$.

\begin{definition}\label{DefPrealgebra}
Let ${(\DD, \tensor, \unit)}$ be a strict monoidal category. A {\em separable semialgebra} is a bi-semigroup
${(D, \nabla, \Delta)}$ in $\DD$
 such that:
\begin{enumerate}
\item (Separable) ${\Delta ; \nabla = id: D \rightarrow D}$
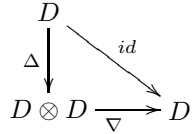
$$\xymatrix{
D \ar[d]_-{\Delta} \ar[rd]^-{id} \\
D\tensor D \ar[r]_{\nabla} & D
}$$
\item (Frobenius) ${(\Delta \tensor id_D); (id_D \tensor \nabla) = \nabla ; \Delta = (id_D \tensor \Delta) ; (\nabla \tensor id_D)}$
$$\xymatrix{
D\tensor D \ar[d]_-{\Delta\tensor id} \ar[r]^-{\nabla} & D \ar[d]^-{\Delta} && 
  D\tensor D \ar[d]_-{id \tensor \Delta} \ar[r]^-{\nabla} & D \ar[d]^-{\Delta} \\
D\tensor D \tensor D \ar[r]_{id \tensor \nabla} & D \tensor D               &&
  D\tensor D \tensor D \ar[r]_{\nabla \tensor id} & D \tensor D
}$$
\end{enumerate}
\end{definition}

Graphically, separability can be expressed as saying that the following two diagrams 
$$
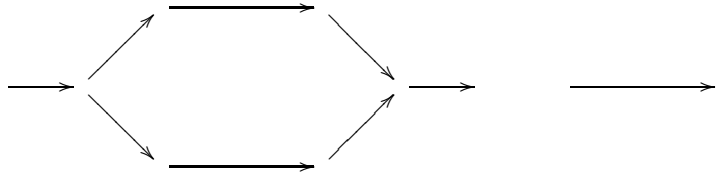
\xymatrix{
       &                & \ar[rr] &  & \ar[rd] &        & &          & & \\
\ar[r] & \ar[ru]\ar[rd] &         &  &         & \ar[r] & & \ar[rr]  & & \\
       &                & \ar[rr] &  & \ar[ru] &        & &          & & \\
}$$
are equal. The authors have found it useful to think of separability as allowing to  {\em pop} the `bubble' on the left.

On the other hand, Frobenius says that the two diagrams below 
$$
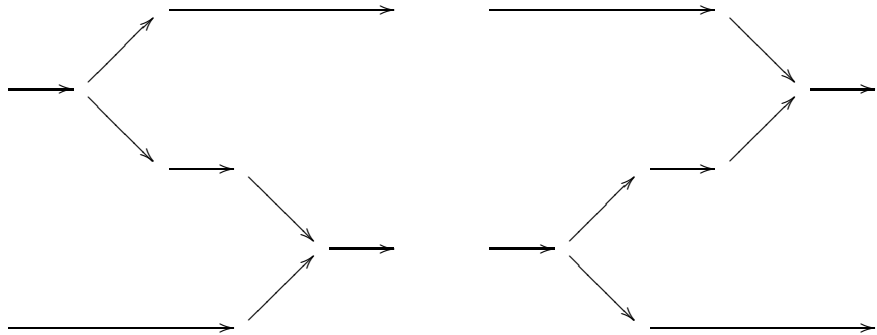
\xymatrix{
       &                & \ar[rrr] &       &        &   & \ar[rrr] &              &        & \ar[rd] &        &  \\
\ar[r] & \ar[ru]\ar[rd] &        &         &        &   &       &                &        &          & \ar[r] &   \\
       &                & \ar[r] & \ar[rd] &        &   &        &                & \ar[r] & \ar[ru] &        & \\
       &                &        &         & \ar[r] &   & \ar[r] & \ar[ru]\ar[rd] &        &         &        & \\ 
\ar[rrr] &              &        & \ar[ru] &        &   &        &                & \ar[rrr] &       &        & \\
}$$
which represent ${(\Delta \tensor id_D); (id_D \tensor \nabla)}$ and ${(id_D \tensor \Delta) ; (\nabla \tensor id_D)}$ respectively, are equal to 
$$
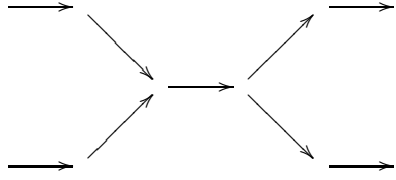
\xymatrix{
\ar[r] & \ar[rd] &         &                & \ar[r] & \\
       &         & \ar[r]  & \ar[ru]\ar[rd] &        & \\ 
\ar[r] & \ar[ru] &         &                & \ar[r] & \\
}$$
which represents ${\nabla ; \Delta}$.

\begin{lemma} If we denote the cospan ${(id:1\rightarrow 1 \leftarrow 1 + 1:\nabla)}$ by ${\Delta:1 \rightarrow 1 + 1}$ then ${(1 + 1, \nabla, \Delta)}$ is a separable semialgebra in $\cospanZ{\sLin}$.
\end{lemma}
\begin{proof}
This is a simple exercise left for the reader. But it is important to mention now that this separable semi-algebra plays an important role in everything that follows.
\end{proof}

\section{A universal property of ${(1 + 1, \nabla, \Delta)}$}

The universal property we discuss in this section was independently observed by Lack on the one hand \cite{Lack04} and by 
Rosebrugh, Sabadini and Walters on the other \cite{RoseSabaWalters05}.

It is important to recall (see Lemma in Chapter~{VII}.5 of \cite{maclane}) the fact that every surjection in $\sLin$ can be factored in a unique way as a composition (satisfying certain conditions) of maps $(id + \nabla + id)$.
(The conditions ensuring uniqueness will not be relevant for us here.)

Let ${F_0:\CC\rightarrow \DD}$ and ${F_1:\opCat{\CC}\rightarrow \DD}$ be functors agreeing on objects. 
We say that ${F_0}$ and ${F_1}$ {\em indulge} a commutaive square 
$$\xymatrix{
\ar[d]_-{\alpha} \ar[r]^-{\beta} & \ar[d]^-{p_1} \\
\ar[r]_-{p_0} &
}$$
if, just as in Definition~\ref{DefCompatible}, ${(F_1 \alpha) ; (F_0 \beta) = (F_0 p_0) ; (F_1 p_1)}$.
(Notice that there is a handedness in this notion. The fact that the functors indulge the square above does not seem to imply that it indulges the square obtained by flipping the same square along its diagonal. That is 
${(F_1 \alpha) ; (F_0 \beta) = (F_0 p_0) ; (F_1 p_1)}$ does not seem to imply
${(F_1 \beta) ; (F_0 \alpha) = (F_0 p_1) ; (F_1 p_0)}$.)

\begin{lemma}\label{LemModular1}
If $F_0$ and $F_1$ indulge the two squares below separately
$$\xymatrix{
\ar[d]_-\alpha \ar[r]^-{\beta} & \ar[d]^-{\alpha'} \ar[r]^-{\beta'} & \ar[d]^-{p'_1} \\
\ar[r]_-{p_0}                  & \ar[r]_-{p'_0} & 
}$$
 then they indulge the rectangle. 
\end{lemma}
\begin{proof}
Trivial.
\end{proof}

There is also a `vertical' version which we shall use when necessary.

\begin{lemma}[See \cite{Lack04}]\label{LemLack}
Let ${F_0:\sLin\rightarrow \DD}$ and ${F_1:\opCat{\sLin}\rightarrow \DD}$ be monoidal functors agreeing on objects. 
Then they are compatible if and only if they indulge the following pushout squares
$$\xymatrix{
1 + 1 \ar[d]_-{\nabla} \ar[r]^-{\nabla} & 1 \ar[d]^-{id} & 
  1 + 1 + 1 \ar[d]_-{\nabla + 1} \ar[r]^-{1 + \nabla} & 1+ 1 \ar[d]^-{\nabla}  &
    1 + 1 + 1 \ar[r]^-{\nabla + 1} \ar[d]_-{1 + \nabla} & 1+ 1 \ar[d]^-{\nabla} \\
1 \ar[r]_-{id}                          & 1              & 
  1 + 1 \ar[r]_-{\nabla} & 1  &
  1 + 1 \ar[r]_-{\nabla} & 1
}
$$
\end{lemma}
\begin{proof}
One direction is trivial.
Consider a pushout of the form below
$$\xymatrix{
k \ar[d]_-{f} \ar[r]^-{g} & n \ar[d]^-{p_1} \\
m \ar[r]_-{p_0} & t
}$$
If $f$ or $g$ are identities then the square is trivially indulged.
So we can assume that $f$ and $g$ are non-trivial compositions. 
Say, ${f = (l_0 + \nabla + l_1); f'}$ and ${g = (l'_0 + \nabla + l'_1); g'}$.
The idea of the proof is to split the pushout into four smaller pushouts as below.
$$\xymatrix{
k \ar[d]_-{l_0 + \nabla + l_1} \ar[rr]^-{l'_0 + \nabla + l'_1} &&   l'_0 + 1 + l'_1 \ar[d] \ar[r]^-{g'} & n \ar[d] \\
l_0 + 1 + l_1 \ar[d]_-{f'} \ar[rr]                             &&    \ar[d] \ar[r] & \ar[d] \\
m \ar[rr]_-{} & & \ar[r] & t
}$$
The inductive hypothesis can deal with two bottom squares and the top right one.
If we can prove that the top left one is indulged then Lemma~\ref{LemModular1} implies that the big pushout is indulged.

So, concerning the top left pushout, the following things can happen:
\begin{enumerate}
\item ${l_0 + 2 \leq l'_0}$, that is, $f$'s first $\nabla$ is strictly to the left of $g$'s,
\item ${l_0 + 1 = l'_0}$, that is, $f$'s first $\nabla$ ``touches" $g$'s but $f$ and $g$ do not start in the same way,
\item ${l_0 = l'_0}$, that is, $f$ and $g$ start in the same way,
\item ${l_0 = l'_0 + 1}$, analogous to  the first item but to the right,
\item ${l_0 \geq l'_0 + 2}$, analogous to the second item.
\end{enumerate}

Consider the first case. 
Let ${k = k_0 + 2 + k_1 + 2 + k_2}$, ${f = (k_0 + \nabla + k'_1);f'}$ and ${g = (k'_0 + \nabla + k_2);g'}$
where ${k'_0 = k_0 + 2 + k_1}$ and ${k'_1 = k_1 + 2 + k_2}$. Then the pushout is calculated as below:
$$\xymatrix{
k_0 + 2 + k_1 + 2 + k_2 \ar[d]_-{k_0 + \nabla + k'_1} \ar[rrr]^-{k'_0 + \nabla + k_2} &&& 
  k_0 +2+ k_1 + 1 + k_2 \ar[d]^-{k_0 + \nabla + k_1 + 1 + k_2}  \\
k_0 + 1 + k_1 + 2 + k_2 \ar[rrr]_-{k_0 + 1 + k_1 + \nabla + k_2} &&&   
  k_0 +2+ k_1 + 1 + k_2   \\
}$$
and it is indulged because it is the sum of trivial pushouts (that are indulged) and moreover $F_0$ and $F_1$ are monoidal so the tensor of indulged squares is indulged. (Should we state a Lemma analogous to Lemma~\ref{LemModular1} but for tensoring squares?.)

For the second case let ${k = k_0 + 1 + 1 + 1 + k_1}$, ${f = (k_0 + \nabla + 1 + k_1); f'}$ and 
${g = (k_0 + 1 + \nabla + k_1); g'}$. In this case the pushout in question is calculated as follows
$$\xymatrix{
k_0 + 1 + 1 + 1 + k_1 \ar[d]_-{k_0 + \nabla + 1 + k_1} \ar[rr]^-{k_0 +1 +\nabla + k_1} && 
  k_0 + 1 + 1 + k_1 \ar[d]^-{k_0 + \nabla + k_1} \\
k_0 + 1 + 1 + k_1 \ar[rr]_-{k_0 + \nabla + k_1} & & k_0 + 1 + k_1 
}$$
Again, the pushout is a sum of two squares that are trivially indulged and one that is indulged by assumption.

   To deal with the third case let ${k = k_0 + 1 + 1 + k_1}$, ${f = (k_0 + \nabla  + k_1); f'}$ and 
${g = (k_0 + \nabla + k_1); g'}$. In this case the pushout in question is calculated as follows
$$\xymatrix{
k_0  + 1 + 1 + k_1 \ar[d]_-{k_0 + \nabla + k_1} \ar[rr]^-{k_0 +\nabla + k_1} && 
  k_0 + 1 + k_1 \ar[d]^-{id} \\
k_0 + 1  + k_1 \ar[rr]_-{id} & & k_0 + 1 + k_1 
}$$
Again, the pushout is a sum of two squares that are trivially indulged and one that is indulged by assumption.

The remaining two cases are analogous.
\end{proof}

\begin{corollary}[See \cite{Lack04} and \cite{RoseSabaWalters05}]\label{CorUniversalPreAlgebra}
Let ${(\DD, \tensor, \unit)}$ be a strict monoidal category with a separable semialgebra ${(D,\nabla,\Delta)}$.
Then there exists a unique strict monoidal functor 
\[ (\cospanZ{\sLin}, +, \ord{0}) \rightarrow (\DD, \tensor, \unit) \]
mapping ${(1 + 1, \nabla, \Delta)}$ to ${(D,\nabla,\Delta)}$.
\end{corollary}
\begin{proof}
The semigroup ${(D,\Delta)}$ is essentially the same thing as a strict monoidal functor 
${F_0:\sLin\rightarrow \DD}$ (mapping $\nabla$ to $\nabla$) while the co-semigroup ${(D,\Delta)}$ is essentially the same thing as a strict monoidal ${F_1:\opCat{\sLin}\rightarrow \DD}$ (mapping $\nabla$ to $\Delta$). 
As $F_0$ and $F_1$ are strict monoidal and coincide on ${1}$, they agree on objects. 
So we are left to prove that $F_0$ and $F_1$ are compatible.
By Lemma~\ref{LemLack} it is enough to check that ${F_0}$ and ${F_1}$ indulge three pushout squares.
But notice that indulgence of these squares is equivalent to Separability and Frobenius.
\end{proof}

\subsection{An alternative proof of Corollary~\ref{CorUniversalPreAlgebra}}

Corollary~\ref{CorUniversalPreAlgebra} can be interpreted as saying that the free monoidal category with a separable semi-algebra is $\cospanZ{\sLin}$. In this short section we sketch a `graphical' proof which makes a lot more evident the relation between the result and the calculation of colimits. 

What is the free monoidal category generated by $\nabla$ and $\Delta$ subject to the equations in Definition~\ref{DefPrealgebra}?  First, given only $\nabla$ we can build diagrams of the form 
$$
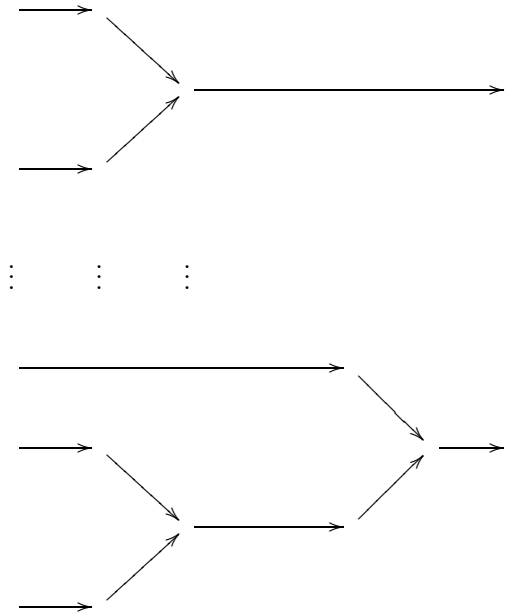
\xymatrix{
\ar[r] & \ar[rd] &        & \\
       &         & \ar[rrrr] & &&  & \\ 
\ar[r] & \ar[ru] &        & &&    \\
\vdots & \vdots  & \vdots \\
\ar[rrrr] &       &        &         & \ar[rd] &      & \\
\ar[r] & \ar[rd] &        &        &         & \ar[r] & \\ 
       &         & \ar[rr] &       & \ar[ru] &        &  \\
\ar[r] & \ar[ru] &  \\
}$$
The associative law says that the order of applying $\nabla$s does not matter so with only $\nabla$s we can can build exactly surjective monotone functions.
Similarly, using only $\Delta$ we can produce exactly the reverses of monotone surjections.
So using both we can produce cospans of monotone surjections. But perhaps we can produce  more? The answer is {\em no}.
If in an expression of $\nabla$s and $\Delta$s a $\Delta$ occurs to the left of a $\nabla$ then only 4 cases can occur.
The first one is when the $\Delta$ and the $\nabla$ do not interact:
$$
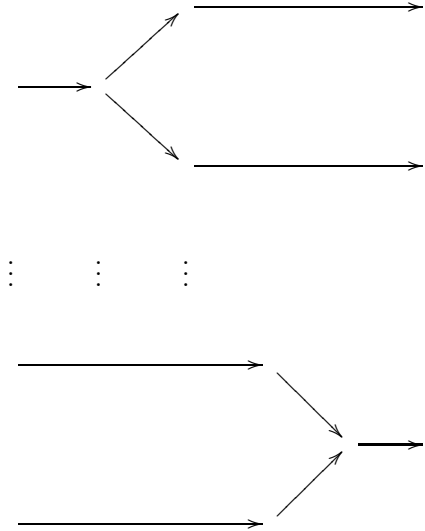
\xymatrix{
       &                & \ar[rrr] & & & \\
\ar[r] & \ar[ru]\ar[rd] &          & & & \\ 
       &                & \ar[rrr] & & & \\
\vdots & \vdots  & \vdots \\
\ar[rrr] &               &         & \ar[rd] &      & \\
          &               &        &          & \ar[r] & \\ 
\ar[rrr] &               &       & \ar[ru] &        &  \\
}$$

The second case is given by the bubble as drawn after Definition~\ref{DefPrealgebra}.
The third and fourth cases are given by the two diagrams representing the expressions in the Frobenius condition and drawn below the bubble after Definition~\ref{DefPrealgebra}.
In all four cases, the $\Delta$s can be moved to the right of the $\nabla$s. 
In the first case trivially, in the second by poping the bubble (separability) and in the third and fourth cases by Frobenius. So the free monoidal category with a separable semi-algebra is $\cospanZ{\sLin}$.

\begin{remark}
It is important to notice that the process of moving $\Delta$s to the right of $\nabla$s is really calculating the pushout involved in the composition of cospans.
\end{remark}

In Section~\ref{SubSecMain} we add 2-dimensional data so that the free monoidal 2-category on this data is $\cospan{\sLin}$.
But first let us extract some more information from Lemma~\ref{LemLack}.

\subsection{Monoidal 2-functors from $\cospan{\sLin}$}

Here we characterize when two functors from $\sLin$ to a 2-category are compatible.
Let us say that a selection of 2-cells $\catc{(\_)}$ indulges a square ${\alpha; p_0 = \beta; p_1}$ if the two equations in Definition~\ref{DefCompatibleFamOf2cells} relating the square and $\catc{(\_)}$ hold.

\begin{lemma}\label{LemIndulgenceFor2cells}
Let ${F_0:\sLin\rightarrow \DD}$ and ${F_1:\opCat{\sLin}\rightarrow \DD}$ be compatible monoidal functors.
Let $\catc{(\_)}$ be a selection of 2-cells satisfying the first two conditions of Definition~\ref{DefCompatibleFamOf2cells}. Then $\catc{(\_)}$ is a compatible selection of 2-cells if and only if it indulges the squares in the statement of Lemma~\ref{LemLack}.
\end{lemma}
\begin{proof}
Analogous to that of Lemma~\ref{LemLack}.
\end{proof}

\section{Adjoint bi-semigroups}

In this section we introduce what we believe are the right lifting to 2-dimensions of the Frobenius and separability conditions.

\begin{definition}
Let ${(\DD, \tensor, \unit)}$ be a strict monoidal 2-category. An {\em adjoint bi-semigroup}
${(X, \nabla, \Delta, \eta, \epsilon)}$ is a bi-semigroup ${(X, \nabla, \Delta)}$ in $\DD_0$ together with 2-cells
${\eta:id_{X\tensor X} \Rightarrow \nabla;\Delta}$ and ${\epsilon:\Delta;\nabla \Rightarrow id_X}$ witnessing that 
${\nabla\dashv \Delta}$.
\end{definition}

We now lift the conditions of separability and Frobenius to the level of adjoint bi-semigroups. We first deal with Frobenius.

\subsection{Frobenius adjoint bi-semigroups}

 In order to justify the definition consider first the following result.

\begin{lemma}\label{LemNablaFrobenius}
Let ${{\bf X} = (X, \nabla, \Delta, \eta, \epsilon)}$ be an adjoint bi-semigroup such that the (1-dimensional) structure ${(X, \nabla, \Delta)}$ satisfies Frobenius as a bi-semigroup. Then the following two items are equivalent:
\begin{enumerate}
\item the mates of the associative laws
$$\xymatrix{
X\tensor X\tensor X \ar[d]_-{\nabla\tensor X} \ar[r]^-{X\tensor\nabla} & X\tensor X \ar[d]^-{\nabla} &&
   X\tensor X\tensor X \ar[d]_-{X\tensor \nabla} \ar[r]^-{\nabla\tensor X} & X\tensor X \ar[d]^-{\nabla} \\
X\tensor X \ar[r]_-{\nabla} & X  && X\tensor X \ar[r]_-{\nabla} & X
}$$
are identity 2-cells

\item ${(\eta\tensor X)\hComp (X\tensor \nabla) = (X\tensor \nabla) \hComp \eta}$ 
and  ${(X\tensor\eta)\hComp(\nabla\tensor X) = (\nabla\tensor X)\hComp\eta}$.
\end{enumerate}
\end{lemma}
\begin{proof}
Consider the mate of one of the associative laws
$$\xymatrix{
X\tensor X \ar[rd]_-{id} \ar[r]^-{\Delta\tensor X} & X\tensor X\tensor X \ar[d]^-{\nabla\tensor X} \ar[r]^-{X\tensor\nabla} & X\tensor X \ar[d]^-{\nabla} \ar[rd]^-{id} & \\
 & X\tensor X \ar[r]_-{\nabla} & X \ar[r]_-{\Delta} & X\tensor X \\
}$$
where inside the triangles we have the 2-cells ${\epsilon\tensor X:(\Delta\tensor X);(\nabla\tensor X) \Rightarrow id}$ and
${\eta:id\Rightarrow \nabla;\Delta}$. Notice that the outside of this diagram is one of the Frobenius laws. 
Now assume that ${(\eta\tensor D)\hComp (D\tensor \nabla) = (D\tensor \nabla) \hComp \eta}$  and calculate:
\[ [(\Delta\tensor X)\hComp (X\tensor \nabla) \hComp\eta] \vComp [(\epsilon\tensor X) \hComp \nabla \hComp \Delta] = \]
\[ = [(\Delta\tensor X)\hComp (\eta\tensor X) \hComp (X\tensor \nabla)]\vComp 
    [(\epsilon\tensor X) \hComp (\Delta\tensor X) \hComp (X\tensor\nabla)]  = \]
\[ = [ ((\Delta\hComp\eta)\tensor X) \hComp (X\tensor\nabla)]\vComp [((\epsilon\hComp\Delta)\tensor X)\hComp (X\tensor\nabla)]  = \]
\[ = [((\Delta\hComp\eta)\tensor X) \vComp ((\epsilon\hComp\Delta)\tensor X)]\hComp (X\tensor\nabla)  = \]
\[ = [ ((\Delta\hComp\eta)\vComp (\epsilon\hComp\Delta)) \tensor X] \hComp (X\tensor\nabla) = 
(\Delta\tensor X) \hComp (X\tensor\nabla) \]
which shows that the mate is the identity 2-cell. Similarly if one assumes that the other equation holds then the corresponding mate is the identity. 

Conversely, assume that the mates of associativity are identity 2-cells and contemplate  the following diagram:
$$\xymatrix{
X^3 \ar[d]_-{X\tensor\nabla} \ar[rd]^-{id} & \\
X^2 \ar[rd]_-{id} \ar[r]^-{X\tensor\Delta} & X^3 \ar[d]^-{X\tensor\nabla} \ar[r]^{\nabla\tensor X} & X^2 \ar[d]^-{\nabla} \ar[rd]^-{id} & \\
                                           & X^2 \ar[r]_-{\nabla}                                  & X \ar[r]_-{\Delta}& X^2\\
}$$
where the triangles are filled with the 2-cells ${X\tensor\eta}$, ${X\tensor \epsilon}$ and ${\eta}$.
Pasting 2-cells one obtains that ${(X\tensor\eta)\hComp (\nabla\tensor X) = (\nabla\tensor X)\hComp\eta }$.
Indeed, one can calculate:
\[ (X\tensor\eta)\hComp (\nabla\tensor X) = \]
\[ = [(X\tensor\eta)\hComp (\nabla\tensor X)]\vComp [(X\tensor\nabla)\hComp (X\tensor\Delta) \hComp (\nabla\tensor X)\hComp\eta]
\vComp [(X\tensor\nabla)\hComp (X\tensor\epsilon)\hComp\nabla\hComp\Delta] = \]
\[ = [(\nabla\tensor X)\hComp \eta]\vComp [(X\tensor\eta)\hComp(\nabla\tensor X)\hComp\nabla\hComp\Delta]\vComp [(X\tensor\nabla)\hComp (X\tensor\epsilon)\hComp\nabla\hComp\Delta] = \]
\[ = [(\nabla\tensor X)\hComp \eta]\vComp [(X\tensor\eta)\hComp(X\tensor \nabla)\hComp\nabla\hComp\Delta]\vComp [(X\tensor (\nabla \hComp \epsilon))\hComp\nabla\hComp\Delta] = \]
\[ = [(\nabla\tensor X)\hComp \eta]\vComp [(X\tensor(\eta\hComp \nabla))\hComp\nabla\hComp\Delta]\vComp [(X\tensor (\nabla \hComp \epsilon))\hComp\nabla\hComp\Delta] = \]
\[ = [(\nabla\tensor X)\hComp \eta]\vComp [[(X\tensor(\eta\hComp \nabla))\vComp (X\tensor (\nabla \hComp \epsilon))]\hComp (\nabla\hComp\Delta) ] =  \]
\[ = (\nabla\tensor X)\hComp \eta \]
The proof of the other equation is analogous.
\end{proof}

Because of this, we find it natural to introduce the following definition.

\begin{definition}\label{DefNablaFrobenius}
Let ${{\bf X} = (X, \nabla, \Delta, \eta, \epsilon)}$ be an adjoint bi-semigroup such that the 1-dimensional structure ${(X, \nabla, \Delta)}$ satisfies Frobenius as a bi-semigroup. We say that {\bf X} satisfies {\em $\nabla$-Frobenius} if the equivalent conditions of Lemma~\ref{LemNablaFrobenius} hold.
\end{definition}

It is interesting and useful to notice that the equalities in the second item of Lemma~\ref{LemNablaFrobenius}  can be thought of as rewrite rules. Indeed, notice that in the notation we have used for expressions with $\Delta$s and $\nabla$s, the left 2-cell of the first equation of item~2 has domain the left hand diagram below
$$
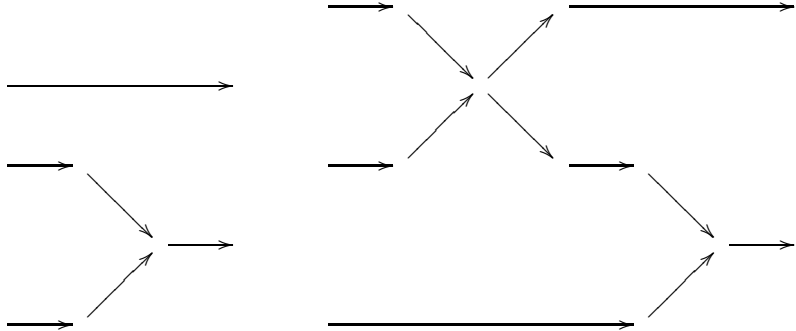
\xymatrix{
         &         &        & & \ar[r] & \ar[rd] &                & \ar[rrr] &         &        &\\
\ar[rrr] &         &        & &        &         & \ar[ru]\ar[rd] &        &         &        & \\
\ar[r]   & \ar[rd] &        & & \ar[r] & \ar[ru] &                & \ar[r] & \ar[rd] &        & \\
         &         & \ar[r] & &        &         &                &        &         & \ar[r] & \\
\ar[r]   & \ar[ru] &        & & \ar[rrrr] &      &                &        & \ar[ru] & \\
}$$
and codomain the right hand diagram below. In other words, the 2-cell {\em pinches} the first two strings. The reader is invited to draw the other 2-cells and exercise in applying the pinching and poping rules.

Back to the lifting of the Frobenius condition, it must be mentioned that  one can prove the following in a way analogous to Lemma~\ref{LemNablaFrobenius}.

\begin{lemma}\label{LemDeltaFrobenius}
Let ${{\bf X} = (X, \nabla, \Delta, \eta, \epsilon)}$ be an adjoint bi-semigroup such that the (1-dimensional) structure ${(X, \nabla, \Delta)}$ satisfies Frobenius as a bi-semigroup. Then the following two items are equivalent:

\begin{enumerate}
\item the mates of 
$$\xymatrix{
X \ar[d]_-{\Delta} \ar[r]^-{\Delta} & X\tensor X \ar[d]^-{\Delta\tensor X} &&
  X \ar[d]_-{\Delta} \ar[r]^-{\Delta} & X\tensor X \ar[d]^-{X\tensor\Delta} \\
X\tensor X \ar[r]_-{X\tensor\Delta} & X\tensor X\tensor X  && X\tensor X \ar[r]_-{\Delta\tensor X} & X\tensor X\tensor X
}$$
are identity 2-cells

\item ${(\Delta\tensor D)\hComp (D\tensor\eta) = \eta \hComp (\Delta\tensor D)}$
and ${(D\tensor\Delta)\hComp(\eta\tensor D) = \eta\hComp (D\tensor\Delta)}$.
\end{enumerate}
\end{lemma}

So, just as in Definition~\ref{DefNablaFrobenius} we say that ${\bf X}$ {\em satisfies $\Delta$-Frobenius} if the equivalent conditions of Lemma~\ref{LemDeltaFrobenius} hold.

\begin{definition}\label{DefFrobeniusAdjBiSemiGrp}
Let ${{\bf X} = (X, \nabla, \Delta, \eta, \epsilon)}$ be an adjoint bi-semigroup such that ${(X, \nabla, \Delta)}$ satisfies Frobenius as a bi-semigroup. We say that {\bf X} satisfies {\em Frobenius} if it satisfies both $\nabla$-Frobenius and $\Delta$-Frobenius.
\end{definition}

\subsection{Separable adjoint bi-semigroups}

In this section we introduce the notion of separable adjoint bi-semigroup and show that for, these semi-groups, 1-dimensional Frobenius implies 2-dimensional Frobenius.

\begin{definition}\label{DefTwoDimSeparability} We say that an adjoint bi-semigroup ${(X, \nabla, \Delta, \eta, \epsilon)}$ is {\em separable} if ${(X, \nabla, \Delta)}$ is separable as a bi-semigroup and moreover ${\epsilon = \iota_{id_X}}$.
\end{definition}

   Notice that in a separable bi-semigroup, ${\eta\hComp \nabla = \iota_{\nabla}}$ and 
${\Delta\hComp\eta = \iota_{\Delta}}$.

\begin{lemma}\label{LemSepImpliesFrob} 
Let ${{\bf X} = (X, \nabla, \Delta, \eta, \epsilon)}$ be an adjoint bi-semigroup such that the 1-dimensional structure
${(X, \nabla, \Delta)}$ satisfies Frobenius. If ${\bf X}$ is separable then ${\bf X}$ satisfies Frobenius.
\end{lemma}
\begin{proof}
Under separability, the triangular identities witnessing that ${\nabla\dashv\Delta}$ become 
${\eta\hComp \nabla = \iota_{\nabla}}$ and ${\Delta\hComp\eta = \iota_{\Delta}}$. To prove $\nabla$-Frobenius we need to show that the mates of associativity are identity 2-cells. In particular, we need to show that 
\[ [(\Delta\tensor X) \hComp (X\tensor \nabla) \hComp \eta]\vComp [(\epsilon\tensor X)\hComp \nabla\hComp\Delta ] \]
is  $\iota_{\nabla ; \Delta}$. Under separability, we need only prove that 
\[ [(\Delta\tensor X) \hComp (X\tensor \nabla) \hComp \eta]\vComp \iota_{\nabla;\Delta} = (\Delta\tensor X) \hComp (X\tensor \nabla) \hComp \eta  \]
is the identity 2-cell $\iota_{\nabla ; \Delta}$. As ${(X, \nabla, \Delta)}$ satisfies Frobenius and 
${\Delta\hComp\eta = \iota_{\Delta}}$, we can calculate:
\[ (\Delta\tensor X) \hComp (X\tensor \nabla) \hComp \eta = \nabla\hComp\Delta\hComp \eta = \nabla\hComp\Delta \]
so, indeed, the mate of associativity is the identity. The other condition is dealt with in an analogous way so 
${\bf X}$ satisfies $\nabla$-Frobenius. To prove $\Delta$-Frobenius one uses the same idea. For example, one of the conditions is proved as follows:
\[ [\eta \hComp (X\tensor\Delta) \hComp (\nabla\tensor X)] \vComp [\nabla\hComp \Delta \hComp (\epsilon \tensor X) ]
= \eta \hComp \nabla \hComp \Delta = \nabla\hComp \Delta \]
so, altogether, ${\bf X}$ satisfies Frobenius.
\end{proof}

   Since the counit is the identity, separable adjoint bi-semigroups will usually be denoted by 
${(X,\nabla, \Delta, \eta)}$.

\subsection{A universal property of $\cospan{\sLin}$}
\label{SubSecMain}

In this section we prove a universal property of $\cospan{\sLin}$ as a monoidal 2-category. For brevity let us introduce the following definition.

\begin{definition}\label{DefCOMOalgebra}
Let ${(\DD, \tensor, \unit)}$ be a strict monoidal  2-category. A {\em Como-algebra} is a separable adjoint bi-semigroup
${(D,\nabla,\Delta, \eta)}$ such that ${(D,\nabla,\Delta)}$ satisfies Frobenius.
\end{definition}

   Alternatively, one can say that a Como-algebra is a structure ${(D,\nabla,\Delta, \eta)}$ such that 
${(D,\nabla,\Delta)}$ is a separable semi-algebra in $\DD_0$ and 
${\eta:id_{D\tensor D} \Rightarrow \nabla;\Delta}$ is a 2-cell satisfying 
\[ \eta\hComp \nabla = \iota_{\nabla} \mbox{\ and\ } \Delta\hComp\eta = \iota_{\Delta} \]
 (essentially saying $\nabla\dashv \Delta$).

By Lemma~\ref{LemSepImpliesFrob} every Como-algebra satisfies Frobenius. Notice also that $1+1$ has an obvious Como-algebra structure: just take ${\eta = \overline{\nabla}}$.

\begin{proposition}\label{PropMainProp}
Let ${(\DD, \tensor, \unit)}$ be a strict monoidal  2-category with a Como-algebra ${(D,\nabla,\Delta, \eta)}$.
Then there exists a unique strict monoidal   2-functor 
\[ (\cospan{\sLin}, +, \ord{0}) \rightarrow (\DD, \tensor, \unit) \]
mapping ${(1 + 1, \nabla, \Delta, \eta)}$ to ${(D,\nabla,\Delta, \eta)}$.
\end{proposition}
\begin{proof}
By Corollary~\ref{CorUniversalPreAlgebra} be have a strict monoidal functor 
${\cospanZ{\sLin} \rightarrow \DD}$ mapping the universal separable semialgebra to the one in $\DD$.
In order to extend this functor to a strict monoidal 2-functor we need a compatible selection of 2-cells satisfying the conditions of Lemma~\ref{LemMonoidalInduced}. That is, a compatible selection $\catc{(\_)}$ satisfying 
${\catc{f + g} = \catc{f} \tensor \catc{g}}$. Now, Definition~\ref{DefCompatibleFamOf2cells} forces $\catc{(\_)}$ on identities and composition. As every map in $\sLin$ is built from $\nabla$ and using tensor and composition, a selection of 2-cells as the one we need is determined by its value ${\catc{\nabla}:id_{D\tensor D} \Rightarrow \nabla ; \Delta}$.
Let us call this selection $\eta$. When does the selection of such a 2-cell induces a compatible selection?
The answer is given by Lemma~\ref{LemIndulgenceFor2cells}. But indulgence of the three distinguished pushouts is equivalent to the validity of the following equations:
\begin{enumerate}
\item ${\eta\hComp \nabla = \iota_{\nabla}}$ and ${\Delta\hComp\eta = \iota_{\Delta}}$ (essentially saying $\nabla\dashv \Delta$)
\item ${(\eta\tensor D)\hComp (D\tensor \nabla) = (D\tensor \nabla) \hComp \eta}$ 
and ${(\Delta\tensor D)\hComp (D\tensor\eta) = \eta \hComp (\Delta\tensor D)}$
\item ${(D\tensor\eta)\hComp(\nabla\tensor D) = (\nabla\tensor D)\hComp\eta}$ 
and ${(D\tensor\Delta)\hComp(\eta\tensor D) = \eta\hComp (D\tensor\Delta)}$
\end{enumerate}
   The first item is exactly separability while the other two items are exactly Frobenius (Definition~\ref{DefFrobeniusAdjBiSemiGrp}).  But a Como-algebra is separable by definition and it always satisfies Frobenius by Lemma~\ref{LemSepImpliesFrob}. So the result follows.
\end{proof}

\section{Como-units}

   Let $\iLin$ be the full subcategory of \Lin\ determined by injective monotone functions. 
The monoidal structure ${(\Lin, +, 0)}$ restricts to \iLin\ and the inclusion ${\iLin\rightarrow\Lin}$ is strict monoidal.
By results in \cite{maclane}, all maps in $\iLin$ are built out of ${!:0\rightarrow 1}$.

\begin{definition}
A {\em unit} in a monoidal category ${(\DD, \tensor, \unit)}$ is an object $X$ in $\DD$ equipped with a map
${u:\unit\rightarrow X}$.
\end{definition}

   The object $1$ in $\iLin$ together with ${!:0\rightarrow 1}$ is the universal object with unit.

\begin{lemma}
The category \iLin\ has strict pullbacks and they interact with $+$.
\end{lemma}

\begin{lemma}\label{LemCompositionalDescriptionOfPb}
Every pullback in $\iLin$ is a composition of trivial pullbacks and pullbacks of the form 
$$\xymatrix{
0 \ar[d]_-{id} \ar[r]^-{id} & 0 \ar[d]^-{!} \\
0 \ar[r]_-{!} & 1 
}$$
\end{lemma}
\begin{proof}
Similar to the proof of Lemma~\ref{LemLack}.
\end{proof}

\begin{definition}
Let ${(\DD, \tensor, \unit)}$ be a monoidal category. A {\em split-unit} is a structure 
${(X, s:\unit\rightarrow X, r:X\rightarrow \unit)}$ such that 
${(X, s:\unit\rightarrow X)}$ is a unit and  ${r:X\rightarrow \unit}$ is such that ${s;r = id_{\unit}}$.
\end{definition}

   The object $1$ has a unique split-unit structure in the monoidal category ${(\cospanZ{\opCat{\iLin}}, +, 0)}$.
Let us denote it by ${(1,!:0\rightarrow 1, ?:1\rightarrow 0)}$.

\begin{corollary}\label{CorSplitPointed}
For every strict monoidal category ${(\DD, \tensor, \unit)}$ and every split-unit
${(X, s:\unit\rightarrow X, r:X \rightarrow \unit)}$ in it, there exists a unique strict monoidal functor ${\cospanZ{\opCat{\iLin}}\rightarrow \DD}$ mapping 
${(1,!:0\rightarrow 1, ?:1\rightarrow 0)}$ to ${(X, s, r)}$.
\end{corollary}
\begin{proof}
The map ${r:X \rightarrow \unit}$ induces a functor strict monoidal ${F_0:\opCat{\iLin}\rightarrow\DD}$ while the map 
${s:\unit\rightarrow X}$ induces a strict monoidal ${F_1:\iLin\rightarrow \DD}$.
The functors clearly agree on objects.
By Lemma~\ref{LemCompositionalDescriptionOfPb}, the functors are compatible if and only if they indulge the pushout
$$\xymatrix{
1 \ar[d]_-{\opCat{!}}\ar[r]^-{\opCat{!}} & 0 \ar[d]^-{\opCat{id}} \\
0 \ar[r]_-{\opCat{id}} & 0
}$$
in $\opCat{\iLin}$. This means exactly that ${s;r = id}$.
\end{proof}

\begin{definition}\label{DefComunit}
Let ${(\DD, \tensor, \unit)}$ be a monoidal 2-category. A {\em Como-unit} is a split-unit
${(X, s:\unit\rightarrow X, l:X \rightarrow \unit)}$ together with a 2-cell 
${\eta:\iota_X \Rightarrow l;s}$ such that ${l\dashv s}$ with unit $\eta$ and counit ${\iota_{id_{\unit}}}$.
\end{definition}

   The split-unit ${(1,!:0\rightarrow 1, ?:1\rightarrow 0)}$ is a Como-unit when considered as an object in  ${(\cospan{\opCat{\iLin}}, +, 0)}$. We denote the unit of the adjunction ${?\dashv !}$ by $\eta$.


In a way analogous to Proposition~\ref{PropMainProp} we obtain the following corollary.

\begin{corollary}\label{CorComunit}
For every strict monoidal 2-category ${(\DD, \tensor, \unit)}$ and Como-unit object ${(X, s, l, \eta)}$ in it, there exists a unique strict monoidal 2-functor
\[ (\cospan{\opCat{\iLin}}, +, 0) \rightarrow (\DD, \tensor, \unit) \]
mapping ${(1,!:0\rightarrow 1, ?:1\rightarrow 0, \eta)}$ to ${(X, s, l, \eta)}$.
\end{corollary}
\begin{proof}
By Corollary~\ref{CorSplitPointed} we have a unique strict monoidal functor 
\[(\cospanZ{\opCat{\iLin}}, +, 0)\rightarrow (\DD_0, \tensor, \unit)\]
mapping the split-unit ${(1,!:0\rightarrow 1, ?:1\rightarrow 0)}$ to ${(X, s, l)}$. 
In order to extend this functor to a strict monoidal 2-functor we need a selection $\catc{(\_)}$ of 2-cells satisfying ${\catc{f+g} = \catc{f}\tensor\catc{g}}$. Such a selection of 2-cells is determined by its value ${\catc{(\opCat{!}:1\rightarrow 0)}:\iota_{X} \Rightarrow l ; s}$.
Naturally, we define ${\catc{(\opCat{!}:1\rightarrow 0)} = \eta}$. Is the resulting selection compatible?
We need to check that $\catc{(\_)}$ indulges all pushout squares in $\opCat{\iLin}$. 
By Lemma~\ref{LemCompositionalDescriptionOfPb} we need only check that it indulges the square in the statement of that lemma. But this says exactly that ${\eta\hComp l = l}$ and ${s\hComp \eta = s}$.
Which means that ${l\dashv s}$ with unit $\eta$ and counit ${\iota_{\unit}}$.
\end{proof}


\end{document}